\documentclass[12pt]{article}
\usepackage[T1]{fontenc}
\usepackage{lmodern}
\usepackage{amsmath,amssymb,amsthm,mathrsfs}
\usepackage{mathtools,xcolor,enumerate}
\usepackage{tikz}
\usepackage{float}
\usetikzlibrary{arrows.meta}
\usepackage[colorlinks=true,allcolors=black]{hyperref}
\usepackage[expansion=false]{microtype}

\newcommand{\R}{\mathbb R}
\newcommand{\tr}{\operatorname{tr}}
\newtheorem{theorem}{Theorem}
\newtheorem{proposition}{Proposition}[section]
\newtheorem{lemma}[proposition]{Lemma}

\theoremstyle{remark}

\numberwithin{equation}{section}
\allowdisplaybreaks[1]
\title{On do Carmo's Question: The Properly Embedded Case}
\author{Gongping Niu}
\date{}
\begin{document}
\maketitle

\begin{abstract}
In this paper, we prove that every complete noncompact two-sided properly embedded constant mean curvature hypersurface of finite Morse index in Euclidean space is minimal. The result holds in every dimension. We assume no bounds on curvature or volume growth. This establishes the properly embedded case of the finite-index version of do Carmo's question.

For nonzero mean curvature, we use the radii of inscribed balls to define a positive weight. This weight gives a polynomial bound on weighted volume, while its reciprocal satisfies a differential inequality involving the Jacobi operator. We combine these properties with stability outside a compact set to obtain finite total area. This leads to a contradiction.
\end{abstract}

\section{Introduction}

A classical conjecture of do Carmo asks whether every complete noncompact two-sided stable hypersurface of constant mean curvature (CMC) immersed in Euclidean space must be minimal \cite{doCarmo1989}. The conjecture predicts a strong geometric consequence of the second variation: stability should rule out nonzero constant mean curvature, without any a priori bounds on curvature or volume growth. If the conjecture is true, that would reduce the study of complete noncompact stable CMC hypersurfaces to the minimal setting.

In this paper, we consider the version for finite-index properly embedded hypersurfaces. Finite Morse index allows unstable directions but implies stability outside a compact set. Thus, the substantial additional assumption is proper embeddedness. We establish the following theorem in arbitrary dimension.

\begin{theorem}\label{thm:main}
Let $n\ge2$, and let $\Sigma^n\subset\R^{n+1}$ be a smooth, connected, complete, noncompact, two-sided, properly embedded hypersurface without boundary and with constant mean curvature. If its Morse index is finite, then its mean curvature vanishes.
\end{theorem}

The theorem imposes no requirements on the second fundamental form, the number of ends, or volume growth. The noncompactness assumption is necessary, since round spheres have nonzero mean curvature and finite Morse index. The finite-index assumption is also essential, even within the properly embedded class. For every $n\ge2$, the round cylinder $\mathbb S^{n-1}(r)\times\R\subset\R^{n+1}$ is complete and properly embedded, with nonzero constant mean curvature and infinite Morse index. 
There are also examples with a wide range of finite topologies. Kapouleas \cite{Kapouleas1990} constructed complete noncompact CMC surfaces by gluing, and Breiner--Kapouleas \cite{BreinerKapouleas2014} refined this construction to obtain many embedded examples of finite topology. More recently, Kleene \cite{Kleene2023} constructed complete embedded CMC surfaces with nonzero mean curvature for every prescribed genus $g\ge0$ and number of ends $k\ge4$. All these complete noncompact nonminimal surfaces in $\R^3$ have infinite Morse index. Thus, finite topology cannot replace the finite-index assumption, even for embedded CMC surfaces.

\medskip

The finite-index problem has been studied extensively in low dimensions. In $\R^3$, it has been proved by Frensel \cite{Frensel1996}; Cheng \cite{Cheng2006} established the corresponding result in ambient dimensions 4 and 5 (see also Elbert--Nelli--Rosenberg \cite{ElbertNelliRosenberg2007} for stable hypersurfaces). Chen--Hong--Li \cite{ChenHongLi2025} and Miranda \cite{Miranda2026} proved the case in $\R^6$, and Wang \cite{Wang2026} obtained the result for hypersurface dimensions $2\le n\le6$. These results apply to immersed hypersurfaces without a properness assumption. Our theorem assumes proper embeddedness and imposes no dimension restriction.

For the properly embedded case, Ilias--Nelli--Soret \cite{IliasNelliSoret2016} proved that every complete noncompact CMC hypersurface with bounded second fundamental form and finite Morse index in Euclidean space is minimal, in arbitrary dimension (see also \cite[Theorem~2.5]{NelliPontuale2026} for the weaker assumption of finite $\delta$-index).  
Our theorem removes the boundedness assumption on the second fundamental form. However, proper embeddedness is essential to our proof. Our argument replaces this uniform geometric control with the inscribed radius. Proper embeddedness ensures that the hypersurface separates Euclidean space into two regions. For nonzero mean curvature, we consider the largest tangent balls contained in the region toward which the mean curvature vector points. Their radii $\rho$ define a weighted measure $\rho\,dV$ with polynomial growth in ambient balls (as shown in Proposition~\ref{prop:weighted-volume}). The reciprocal $\mu=\rho^{-1}$ satisfies $\mu\ge H$ and the global distributional inequality $\Delta\mu+|A|^2\mu\ge nH\mu^2$ (Proposition~\ref{prop:mu-distributional}). The potential in this inequality agrees with that of the Jacobi operator, which allows us to combine it with stability outside a compact set. The resulting integral estimates imply finite total area, contrary to the infinite-volume result (\cite[Appendix]{Frensel1996}).

\medskip

The analytic part of our argument also applies in Riemannian ambient manifolds. It requires a weight $\mu\ge c>0$ with a corresponding differential inequality for the full Jacobi operator and polynomial growth of the measure $\mu^{-1}\,dV$ on intrinsic balls. In order to lead to a contradiction, the ambient geometry also needs to ensure that complete noncompact CMC hypersurfaces have infinite volume.

In Euclidean space, the remaining issue is to construct such a weight without proper embeddedness. Once this weight is available, intrinsic distance cutoffs allow the argument to apply to complete immersed hypersurfaces. Constructing such a function for general complete noncompact two-sided CMC immersions with nonzero mean curvature would establish the full finite-index version of do Carmo's conjecture.

\subsection*{Notation and preliminaries}

Throughout the proof, $\Sigma$ satisfies the geometric assumptions of
Theorem~\ref{thm:main}. We identify it with its embedded image. We use
the following notation.
\begin{list}{}{%
 \setlength{\labelwidth}{3.4cm}%
 \setlength{\leftmargin}{3.7cm}%
 \setlength{\labelsep}{0.3cm}%
 \setlength{\itemsep}{0.35em}%
 \setlength{\parsep}{0pt}%
 \renewcommand{\makelabel}[1]{#1\hfil}}
\item[$\Sigma^n$] the hypersurface in $\R^{n+1}$;
\item[$g, dV$] the induced metric and volume measure on $\Sigma$;
\item[$\nabla$] the Levi-Civita connection on $(\Sigma,g)$.
For a function $f$ on $\Sigma$, $\nabla f$ is its gradient;
\item[$\Delta$] the Laplacian on $(\Sigma,g)$;
\item[$\exp_x$] the Riemannian exponential map of $\Sigma$ at $x$;
\item[$N$] a smooth global unit normal;
\item[$A$] the shape operator, with $A(v)=-D_vN$. We also write
$A(v,w)=\langle A(v),w\rangle$;
\item[$H$] the constant mean curvature, with $H=\frac1n\tr A$;
\item[$\kappa_i$] the principal curvatures, ordered as
$\kappa_1\le\cdots\le\kappa_n$;
\item[$B_R(a),B_R$] the ambient ball of radius $R$ centered
at $a$, and the ball $B_R(0)$, respectively;
\item[$\omega_m$]  the volume of the unit ball in $\R^m$.
\end{list}
We omit $dV$ from integrals when the measure is clear. Both
$\operatorname{Area}$ and $\operatorname{Vol}_n$ denote volume with
respect to $dV$.

If $H\ne0$, we choose $N$ so that $H>0$. We keep this orientation
throughout the proof. By the separation theorem for connected properly
embedded hypersurfaces without boundary, we can
choose the complementary domain $\Omega$ into which $N$ points.
Hence $\partial\Omega=\Sigma$.

Note that the principal curvatures satisfy
\[
 |A|^2=\sum_{i=1}^n\kappa_i^2\ge nH^2.
\]
We define the Jacobi quadratic form by
\begin{equation}\label{eq:quadratic}
 Q(f)=\int_\Sigma\bigl(|\nabla f|^2-|A|^2f^2\bigr),
 \qquad f\in C_c^\infty(\Sigma).
\end{equation}
The Morse index is the supremum of the dimensions of subspaces of
$C_c^\infty(\Sigma)$ on which $Q$ is negative definite.
If the Morse index is finite, we can choose a compact set
$E\subset\Sigma$ such that
\begin{equation}\label{eq:exterior-stability}
 \int_\Sigma |A|^2f^2\le\int_\Sigma|\nabla f|^2,
 \qquad f\in C_c^\infty(\Sigma\setminus E).
\end{equation}
Otherwise, we could choose arbitrarily many negative directions with
pairwise disjoint compact supports. Their span would be negative definite, contradicting the finite-index assumption.

By density, \eqref{eq:exterior-stability} also holds for compactly supported Lipschitz functions on $\Sigma\setminus E$.

\medskip
To motivate the weighted argument below, we first give a cutoff proof of minimality under a polynomial volume growth assumption (see also Alencar--do Carmo \cite{AlencarDoCarmo1993,AlencarDoCarmo1995Erratum}). Assume that there exist constants $C,d>0$ such that
\[
V(R):=\operatorname{Area}(\Sigma\cap B_R)\le CR^d \qquad \text{for all }R\ge1.
\]
Suppose that $H\ne0$. Choose $\zeta\in C^\infty(\Sigma)$ with $0\le\zeta\le1$, vanishing on a neighborhood of $E$ and equal to one outside a larger compact set $E'$. Fix a real exponent $p>\max\{d,2\}$, and let $\eta_R$ be the restriction to $\Sigma$ of a smooth ambient cutoff satisfying
\[
0\le\eta_R\le1,\qquad \eta_R=1\text{ on }B_R,\qquad \eta_R=0\text{ outside }B_{2R},\qquad |\nabla\eta_R|\le C/R.
\]
Properness ensures that $\chi_R:=\zeta\eta_R$ has compact support in $\Sigma\setminus E$.
Set $I_R:=\int_\Sigma\chi_R^p\,dV$. Testing \eqref{eq:exterior-stability} with $\chi_R^{p/2}$ gives
\[
\begin{aligned}
nH^2 I_R
&\le \frac{p^2}{4}\int_\Sigma\chi_R^{p-2}|\nabla\chi_R|^2\,dV\\
&\le \frac{p^2}{4}I_R^{1-2/p}\left(\int_\Sigma|\nabla\chi_R|^p\,dV\right)^{2/p}.
\end{aligned}
\]
Consequently,
\[
I_R\le\left(\frac{p^2}{4nH^2}\right)^{p/2}\int_\Sigma|\nabla\chi_R|^p\,dV.
\]
For large $R$, we have $\eta_R=1$ on $E'$ and $\zeta=1$ on the support of $\nabla\eta_R$, so the two terms in $\nabla\chi_R=\eta_R\nabla\zeta+\zeta\nabla\eta_R$ have disjoint supports. The volume growth assumption therefore yields
\[
\begin{aligned}
\int_\Sigma|\nabla\chi_R|^p\,dV
&\le \int_\Sigma|\nabla\zeta|^p\,dV+\frac{C^p}{R^p}V(2R)\\
&\le C_0+C_1R^{d-p}=C_0+o(1),
\end{aligned}
\]
where the constants are independent of $R$. Since $\eta_R\to1$ pointwise, Fatou's lemma gives
\[
\int_\Sigma\zeta^p\,dV<\infty.
\]
As $\zeta=1$ outside $E'$ and $E'$ has finite area, it follows that $\operatorname{Area}(\Sigma)<\infty$. This contradicts Frensel's infinite-volume theorem \cite[Appendix, Theorem~1]{Frensel1996}. Hence $H=0$.

\subsection*{Outline of the proof}

Suppose that $H>0$. We will find a contradiction.
We look for a locally Lipschitz function $\mu\ge H>0$ with the following two properties.
\begin{itemize}
 \item The function $\mu$ satisfies the weighted stability inequality
 \[
  nH\int_\Sigma\mu^3\varphi^2
  \le\int_\Sigma\mu^2|\nabla\varphi|^2
 \]
 for every compactly supported Lipschitz function $\varphi$ on
 $\Sigma\setminus E$ (Proposition~\ref{prop:weighted-stability}).

 \item Its reciprocal $\rho=\mu^{-1}$ satisfies the geometric volume
 estimate
 \begin{equation}\label{eq:intro-volume-estimate}
  \int_{\Sigma\cap B_R}\rho
  \le(n+1)\omega_{n+1}(R+H^{-1})^{n+1}
 \end{equation}
 for every $R>0$ (Proposition~\ref{prop:weighted-volume}).
\end{itemize}
To construct this weight, we consider the inscribed radius and define
\[
 \mu(x)=\sup_{y\in\Sigma\setminus\{x\}}
 \frac{2\langle y-x,N(x)\rangle}{|y-x|^2}.
\]
Geometrically, the reciprocal $\rho$ is the radius of the largest tangent ball on the
side into which $N$ points. Thickening $\Sigma\cap B_R$ along the inward normals by the local distance $\rho(x)$ produces a region without overlap. Comparing its volume with that of $B_{R+H^{-1}}$ gives \eqref{eq:intro-volume-estimate}.

The maximal inscribed ball at a point is limited either by the local curvature or by contact with a second point of $\Sigma$ (see the two cases in the proof of Proposition~\ref{prop:mu-distributional}). Following Brendle's two-point approach \cite{Brendle2015}, we use this geometry and the CMC equation to show that $\mu$ satisfies
\begin{equation}\label{eq:intro-mu}
\Delta\mu+|A|^2\mu\ge nH\mu^2
\end{equation}
in the sense of distributions on all of $\Sigma$.

The potential $|A|^2$ matches that in exterior stability \eqref{eq:exterior-stability}, allowing us to combine the two inequalities to obtain the weighted stability estimate above. The rest of the proof parallels the argument under the polynomial volume growth assumption (with \eqref{eq:intro-volume-estimate} providing the required control of the cutoff error). We again obtain finite total area, contradicting Frensel's infinite-volume theorem.

\section{Geometry of the inscribed radius}\label{sec:geometry}

In this section, we recall the inscribed-radius geometry used in the
two-point arguments of Andrews \cite{Andrews2012} and Brendle \cite{Brendle2015}
and prove a weighted volume estimate.
Note that we do not need the mean curvature to be constant in this
section. We only need $\Sigma$ to be a properly embedded hypersurface
with $H(x)\ge H_0>0$.

We define
\begin{equation*}
 \mu(x)=\sup_{y\in\Sigma\setminus\{x\}}
 \frac{2\langle y-x,N(x)\rangle}{|y-x|^2},
 \qquad \rho(x)=\frac1{\mu(x)}.
\end{equation*}
We take a unit vector $v\in T_x\Sigma$. Along a curve through $x$ with
velocity $v$, we use Taylor expansion to see that the quotient tends to
$A_x(v,v)$. Therefore, we have
\begin{equation}\label{eq:radius-bounds}
\mu(x)\ge\kappa_n(x)\ge H(x)\ge H_0>0,
\qquad
0<\rho(x)\le H(x)^{-1}\le H_0^{-1}.
\end{equation}
By Lemma~\ref{lem:radius-regularity}, we know that $\mu$ is finite and
locally Lipschitz.

We next recall the standard interior-ball geometry.
We include the short proof for completeness.

\begin{lemma}\label{lem:inscribed-balls}
For each $x\in\Sigma$, $\rho(x)$ is the radius of the largest open ball
contained in $\Omega$ and tangent to $\Sigma$ at $x$. This ball is
\[
 B_{\rho(x)}\bigl(x+\rho(x)N(x)\bigr).
\]
If a distinct point $y_0\in\Sigma$ realizes the supremum defining
$\mu(x)$, then $y_0$ lies on the boundary of this ball and satisfies
\begin{equation}\label{eq:normal-reflection}
 N(y_0)=N(x)-\mu(x)(y_0-x).
\end{equation}
In particular, for every $0<t<\rho(x)$, the distance from $x+tN(x)$ to $\Sigma$ is $t$,
and $x$ is the unique nearest point.
\end{lemma}

\begin{proof}
We fix $x\in\Sigma$.
For $y\ne x$ and $0<t\le\rho(x)$, we consider the distance between
$x+tN(x)$ and $y$. By the definition of $\mu(x)$, we have
\begin{align}
 |y-x-tN(x)|^2-t^2\label{eq:nearest-point-eq}
 &=|y-x|^2-2t\langle y-x,N(x)\rangle\\
 &\ge(1-t\mu(x))|y-x|^2.\label{eq:nearest-point}
 \end{align}
Therefore, when $t=\rho(x)$, the open ball is disjoint from $\Sigma$
and lies in $\Omega$.

We now show that $\rho(x)$ is maximal. An open ball contained in
$\Omega$ and tangent to $\Sigma$ at $x$ has center $x+rN(x)$ and radius
$r$.
If $r>\rho(x)$, we can choose a point $y\in\Sigma\setminus\{x\}$ with
$2r\langle y-x,N(x)\rangle>|y-x|^2$. By \eqref{eq:nearest-point-eq}
with $t=r$, we have
$$ |y-x-rN(x)|< r.$$
This point lies inside $B_r(x+rN(x))$, so that ball is not contained in
$\Omega$.

\medskip

If $y_0\ne x$ realizes $\mu(x)$, then equality holds in
\eqref{eq:nearest-point} at $t=\rho(x)$. Thus $y_0$ lies on the boundary
of the ball. Since the ball lies in $\Omega$, it is tangent to $\Sigma$
at $y_0$. Its inward unit normal there is $N(y_0)$. Therefore, we have
\eqref{eq:normal-reflection}.

\medskip

Finally, if $0<t<\rho(x)$, the right-hand side of
\eqref{eq:nearest-point} is strictly positive for every $y\ne x$.
Every such point is therefore farther than $t$ from $x+tN(x)$,
while $x$ is at distance $t$. Therefore, we have proved the last assertion.
\end{proof}

\begin{proposition}\label{prop:weighted-volume}
For every $R>0$,
\begin{equation*}
\int_{\Sigma\cap B_R}\rho\,dV
\le(n+1)\omega_{n+1}(R+H_0^{-1})^{n+1}.
\end{equation*}
\end{proposition}

\begin{proof}
We consider the normal map
\[
 F(x,t)=x+tN(x),\qquad
 \mathcal D=\{(x,t):x\in\Sigma,\ 0<t<\rho(x)\}.
\]
By Lemma~\ref{lem:inscribed-balls}, we know that $F$ is injective on
$\mathcal D$.
We fix $(x,t)\in\mathcal D$ and choose an orthonormal basis
$e_1,\ldots,e_n$ of principal directions at $x$. Since
$D_{e_i}N=-\kappa_i(x)e_i$, we have
\[
 DF_{(x,t)}(e_i,0)
 =(1-t\kappa_i(x))e_i,
 \qquad
 DF_{(x,t)}(0,1)=N(x).
\]
Note that
$\kappa_i(x)\le\mu(x)=\rho(x)^{-1}$ and $0<t<\rho(x)$. Therefore, we have
\[
 1-t\kappa_i(x)\ge 1-\frac{t}{\rho(x)}>0.
\]
Thus the Jacobian of $F$ is
\[
 J_F(x,t)
 =\prod_{i=1}^n(1-t\kappa_i(x))
 \ge \left(1-\frac{t}{\rho(x)}\right)^n.
\]
By integrating in $t$, we have
\[
 \int_0^{\rho(x)}J_F(x,t)\,dt
 \ge
 \int_0^{\rho(x)}
 \left(1-\frac{t}{\rho(x)}\right)^n\,dt
 =\frac{\rho(x)}{n+1}.
\]

We now fix $R>0$, restrict the base point to $\Sigma\cap B_R$, and define
\[
 \mathcal D_R
 =\{(x,t):x\in\Sigma\cap B_R,\ 0<t<\rho(x)\}.
\]
For every $(x,t)\in\mathcal D_R$, we use
\eqref{eq:radius-bounds} to get
\[
 |F(x,t)|
 \le |x|+t
 <R+\rho(x)
 \le R+H_0^{-1}.
\]
Therefore, $F(\mathcal D_R)\subset B_{R+H_0^{-1}}$, and we have
\begin{align*}
 \omega_{n+1}(R+H_0^{-1})^{n+1}
 &\ge \operatorname{Vol}_{n+1}\bigl(F(\mathcal D_R)\bigr)\\
 &=\int_{\mathcal D_R}J_F(x,t)\,dV(x)\,dt\\
 &=\int_{\Sigma\cap B_R}
      \int_0^{\rho(x)}J_F(x,t)\,dt\,dV(x)\\
 &\ge\frac1{n+1}\int_{\Sigma\cap B_R}\rho(x)\,dV(x).
\end{align*}
\end{proof}

\section{The differential inequality for \texorpdfstring{$\mu$}{mu}}\label{sec:differential-inequality}

From now on, we again assume that the mean curvature $H$ is a positive constant. We first prove the
distributional inequality in Proposition~\ref{prop:mu-distributional}.
We then combine it with exterior stability to prove Theorem~\ref{thm:main}.

\medskip

We first study the local regularity of $\mu$. 
Brendle proves local Lipschitz continuity and semiconvexity on
$\{\mu>\kappa_n\}$ (cf. \cite[Proposition~2.1]{Brendle2015}).
To include points where $\mu=\kappa_n$, we use a local form of the
diagonal construction
of Andrews--Langford--McCoy
(cf. \cite[Section~4]{AndrewsLangfordMcCoy2013}).

\begin{lemma}\label{lem:radius-regularity}
Let $\Sigma\subset\mathbb R^{n+1}$ be a smooth properly embedded
hypersurface with a smooth unit normal $N$. The function
\[
 \mu(x)=\sup_{y\in\Sigma\setminus\{x\}}
 \frac{2\langle y-x,N(x)\rangle}{|y-x|^2}
\]
is finite, locally Lipschitz, and locally semiconvex.
\end{lemma}

\begin{proof}
We fix a point of $\Sigma$. After a rigid motion, we choose a graph
parametrization
\[
 X:B_{5r}^n\longrightarrow\Sigma,
 \qquad X(z)=(z,u(z)),
\]
around this point.

We first consider nearby second points. We write them as
$X(z+t\omega)$ for $z\in B_r^n$, $\omega\in S^{n-1}$, and $0<t\le3r$.
We define
\[
 G_{t,\omega}(z)=
 \frac{2\langle X(z+t\omega)-X(z),N(X(z))\rangle}
 {|X(z+t\omega)-X(z)|^2}.
\]
Since $\langle DX(z)\omega,N(X(z))\rangle=0$, we use Taylor's
formula to write
\begin{align*}
 X(z+t\omega)-X(z)
 &=t\int_0^1 DX(z+st\omega)\omega\,ds,\\
 \langle X(z+t\omega)-X(z),N(X(z))\rangle
 &=t^2\int_0^1(1-s)
 \langle D^2X(z+st\omega)[\omega,\omega],N(X(z))\rangle\,ds.
\end{align*}
Therefore, we have
\[
 G_{t,\omega}(z)=
 \frac{2\displaystyle\int_0^1(1-s)
 \langle D^2X(z+st\omega)[\omega,\omega],N(X(z))\rangle\,ds}
 {|X(z+t\omega)-X(z)|^2/t^2}.
\]
Since $X(z)=(z,u(z))$ and $|\omega|=1$, we have
\[
 |X(z+t\omega)-X(z)|\ge t.\]
Therefore, we can extend the function smoothly in $(z,t,\omega)$ to $t=0$, with
\[
 G_{0,\omega}(z)=
 \frac{\langle D^2X(z)[\omega,\omega],N(X(z))\rangle}
 {|DX(z)\omega|^2}.
\]
The value $G_{0,\omega}(z)$ is the normal curvature in the unit direction $DX(z)\omega/|DX(z)\omega|$. By smoothness and compactness, the extended family has uniform $C^2$ bounds in $z$ on $\overline B_r^n$.

We next consider second points outside $X(B_{2r}^n)$. We set
\[
 \mathcal F=\Sigma\setminus X(B_{2r}^n),
 \qquad
 G_y(z)=
 \frac{2\langle y-X(z),N(X(z))\rangle}{|y-X(z)|^2},
 \quad y\in\mathcal F.
\]
By properness, $\mathcal F$ is closed in $\mathbb R^{n+1}$.
It has positive distance from the compact set $X(\overline B_r^n)$. We conclude that $\{G_y: y\in \mathcal{F}\}$ has uniform $C^2$ bounds on $\overline B_r^n$.

For $z\in B_r^n$ and
$y=X(\xi)\in X(B_{2r}^n)\setminus\{X(z)\}$, we take $t\in(0,3r),
 \omega\in S^{n-1}$ and write $y=X(z+t\omega)$. Thus the two families cover all second
points $y\ne X(z)$. 
We therefore obtain
\[
\mu(X(z))
=
\sup\left\{
G_{t,\omega}(z), G_y(z)
:
(t,\omega)\in[0,3r]\times S^{n-1},\;
y\in\mathcal{F}
\right\}.
\]

Both families have uniform $C^2$ bounds. 
By taking the supremum, we conclude that $\mu$ is finite, Lipschitz continuous and semiconvex.

\end{proof}

\begin{proposition}[cf.\ {\cite[Corollary~2.4]{Brendle2015}}]\label{prop:mu-distributional}
The reciprocal inscribed radius $\mu$ satisfies
\begin{equation}\label{eq:mu-distributional}
 \Delta\mu+|A|^2\mu\ge nH\mu^2
\end{equation}
in the sense of distributions on $\Sigma$.
\end{proposition}

\begin{proof}
By Lemma~\ref{lem:radius-regularity} and Alexandrov's theorem
(see\ \cite[Section~6.4]{EvansGariepy1992}), $\mu$ has a second-order
Taylor expansion at almost every point. We fix such a point $x$ and write
\[
 \mu(\exp_x v)=\mu(x)+\langle\nabla\mu(x),v\rangle
 +\frac12 B_x(v,v)+o(|v|^2).
\]
Here $B_x$ is the Hessian of $\mu$ at $x$ in the sense of Alexandrov.
We prove
\begin{equation}\label{eq:alexandrov-inequality}
 \operatorname{tr}_g B_x+|A(x)|^2\mu(x)\ge nH\mu(x)^2.
\end{equation}

\bigskip

\noindent\emph{Case 1: $\mu(x)>\kappa_n(x)$.}

Hence, the supremum defining $\mu(x)$ is attained at a point $y\neq x$. 
Indeed, the
quotient is at most $2/|y-x|$ and $\mu(x)\geq H>0$. A maximizing
sequence therefore stays in a bounded ambient set. By properness, we can
pass to a convergent subsequence. On the other hand, 
the diagonal expansion gives
\[
\limsup_{y\to x}\frac{2\langle y-x,N(x)\rangle}{|y-x|^2}=\kappa_n(x).
\]
Since $\mu(x)>\kappa_n(x)$, the maximizing sequence cannot converge to $x$.

\medskip

Denote $q=y-x$. By Lemma~\ref{lem:inscribed-balls}, 
\[
 \langle q,N(x)\rangle=\frac{\mu|q|^2}{2},
 \qquad
 N(y)=N(x)-\mu q.
\]

To compare second-order terms at the contact pair, we define
\[
 Z(x',y')=
 \frac{\mu(x')}{2}|y'-x'|^2-
 \langle y'-x',N(x')\rangle,
 \qquad x'\text{ near }x,\quad y'\in\Sigma.
\]
Then by definition of $\mu$, $Z\geq0$ and $Z(x,y)=0$. The first-order terms of $Z$ vanish. We choose an orthonormal basis $e_1,\ldots,e_n$ of principal directions
at $x$. By varying $x'$ in the direction $e_i$ while keeping $y'=y$,
we obtain
\[
 0=\frac{|q|^2}{2}\nabla_{e_i}\mu
   -(\mu-\kappa_i)\langle q,e_i\rangle.
\]
Since $\mu>\kappa_n\ge\kappa_i$, we have
\begin{equation}\label{eq:contact-gradient}
 \langle q,e_i\rangle
 =\frac{|q|^2}{2(\mu-\kappa_i)}\nabla_{e_i}\mu.
\end{equation}
Next we denote by $\bar e_i$ the orthogonal reflection of $e_i$ across $q^\perp$.
\[
\bar e_i
=e_i-2\frac{\langle q,e_i\rangle}{|q|^2}q.
\]
The inscribed ball is tangent to $\Sigma$ at both $x$ and $y$.
The same reflection therefore maps $N(x)$ to $N(y)$. Since it preserves inner products, the vectors $\bar e_i$ form an
orthonormal basis of $T_y\Sigma$.

For each $i$, we choose intrinsic geodesics $x_i(t)$ and $y_i(t)$ with initial
data $(x,e_i)$ and $(y,\bar e_i)$. Along these curves, we write
\[
 Z(x_i(t),y_i(t))=\frac12 b_i t^2+o(t^2),
 \qquad b_i\ge0.
\]
We next follow Brendle's two-point argument and give the calculation in our
conventions.
By the reflection formula and the geodesic equations, we have
\begin{align*}
 (y_i-x_i)'(0)
 &=-\frac{2\langle q,e_i\rangle}{|q|^2}\,q,\\
 (y_i-x_i)''(0)
 &=A_y(\bar e_i,\bar e_i)N(y)-\kappa_iN(x).
\end{align*}

Therefore, using the second-order Taylor expansion of $\mu$ at $x$, a direct computation gives
\begin{align*}
 b_i
 ={}&\frac{|q|^2}{2}
 \left(
 B_x(e_i,e_i)+\mu\kappa_i^2-\mu^2\kappa_i
 \right)\\
 &-4(\nabla_{e_i}\mu)\langle q,e_i\rangle
 +\frac{4(\mu-\kappa_i)}{|q|^2}\langle q,e_i\rangle^2\\
 &+\kappa_i-A_y(\bar e_i,\bar e_i)
 +\left\langle q,(\nabla_{e_i}A)(e_i)\right\rangle.
\end{align*}
Using \eqref{eq:contact-gradient}, we get
\begin{align*}
 b_i
 ={}&\frac{|q|^2}{2}
 \left(
 B_x(e_i,e_i)+\mu\kappa_i^2-\mu^2\kappa_i
 -\frac{2(\nabla_{e_i}\mu)^2}{\mu-\kappa_i}
 \right)\\
 &+\kappa_i-A_y(\bar e_i,\bar e_i)
 +\left\langle q,(\nabla_{e_i}A)(e_i)\right\rangle.
\end{align*}

By orthonormality and the Codazzi equations, we obtain
\[
 \sum_i A_y(\bar e_i,\bar e_i)=nH(y),
 \qquad
 \sum_i(\nabla_{e_i}A)(e_i)=n\nabla H(x).
\]
Since $\nabla H=0$ and $H(x)=H(y)=H$, the last three terms cancel
after summing. We obtain
\begin{equation}\label{eq:two-point-second-variation}
 \sum_i b_i
 =\frac{|q|^2}{2}
 \left(
 \operatorname{tr}_g B_x+|A|^2\mu-nH\mu^2
 -2\sum_i\frac{(\nabla_{e_i}\mu)^2}{\mu-\kappa_i}
 \right).
\end{equation}
Since $q\ne0$ and $b_i\ge0$,
\eqref{eq:two-point-second-variation} implies that
\[
 \operatorname{tr}_g B_x+|A|^2\mu
 \ge nH\mu^2
 +2\sum_i\frac{(\nabla_{e_i}\mu)^2}{\mu-\kappa_i}.
\]
Hence we obtain \eqref{eq:alexandrov-inequality} at every point of $\{\mu>\kappa_n\}$ where $\mu$
has a second-order Taylor expansion.

\bigskip

\noindent\emph{Case 2: $\mu(x)=\kappa_n(x)$.}

We choose a unit eigenvector $v$ of $A_x$ with
eigenvalue $\kappa_n(x)=\mu(x)$ and extend it to a smooth local unit tangent vector
field $V$ with $\nabla V(x)=0$ on $\Sigma$. 
Then the smooth function $A(V,V)$ lies below $\mu$ and equals it at $x$:
\[
 A(V,V)\leq\kappa_n\leq\mu
 \quad\text{near }x,
 \qquad A_x(v,v)=\mu(x).
\]
Therefore, at $x$, by comparing their second-order expansions,
\[
 B_x\ge\nabla^2(A(V,V))(x), \qquad  \operatorname{tr}_g B_x\ge\Delta(A(V,V))(x).
\]

By the tensor Simons identity
\cite[Appendix~B]{Simon2018}:
\begin{equation*}
 \Delta A=nH A^2-|A|^2A.
\end{equation*}
At $x$, the terms involving $\nabla V$ vanish:
\[
\Delta(A(V,V))(x)
=(\Delta A)(v,v)+2A(\Delta V,v).
\]
On the other hand, differentiating $|V|^2=1$ gives
$\langle\Delta V,V\rangle=-|\nabla V|^2=0$ at $x$. So the remaining vector-field term also vanishes:
\[
 2A(\Delta V,V)=2\mu(x)\langle\Delta V,V\rangle=0.
\]

Consequently, we obtain
\[
 \operatorname{tr}_g B_x\ge\Delta(A(V,V))(x)
 =nH\mu(x)^2-|A(x)|^2\mu(x).
\]

We have proved \eqref{eq:alexandrov-inequality} at almost every point of
$\Sigma$. Since $\mu$ is locally semiconvex, its distributional
Laplacian is a signed Radon measure. Its absolutely continuous part is
$(\operatorname{tr}_g B_x)\,dV$, and its singular part is nonnegative
(see \cite[Section~6.4]{EvansGariepy1992}). Thus
\[
 \Delta\mu\ge(\operatorname{tr}_g B_x)\,dV.
\]
By combining the two cases together, we obtain \eqref{eq:mu-distributional} in the sense of distributions on $\Sigma$.
\end{proof}

In particular, we obtain
\begin{equation}\label{eq:mu-integrated}
 -\int_\Sigma\langle\nabla\mu,\nabla\psi\rangle\,dV
 +\int_\Sigma |A|^2\mu\psi\,dV
 \ge nH\int_\Sigma\mu^2\psi\,dV
 \qquad(0\le\psi\in C_c^\infty(\Sigma)).
\end{equation}
The positive function $\mu$ satisfies \eqref{eq:mu-distributional}
with the same Jacobi potential $|A|^2$ as exterior stability
\eqref{eq:exterior-stability}. We test the differential inequality with
$\mu\varphi^2$ and stability with $\mu\varphi$ to obtain the following
weighted stability estimate.

\begin{proposition}\label{prop:weighted-stability}
Suppose that \eqref{eq:exterior-stability} holds outside a compact set
$E\subset\Sigma$. Then every compactly supported
Lipschitz function $\varphi$ on $\Sigma\setminus E$ satisfies
\begin{equation}\label{eq:stability-cancellation}
 nH\int_\Sigma\mu^3\varphi^2\,dV
 \le \int_\Sigma\mu^2|\nabla\varphi|^2\,dV.
\end{equation}
Consequently, for every real $p\ge8$ and every compactly supported
Lipschitz function $\chi$ on $\Sigma\setminus E$ with $0\le\chi\le1$,
\begin{equation}\label{eq:arbitrary-power}
 \int_\Sigma\mu^3\chi^p\,dV
 \le
 \left(\frac{p^2}{4n}\right)^{p/2}H^{4-p}
 \int_\Sigma\rho|\nabla\chi|^p\,dV.
\end{equation}
\end{proposition}

\begin{proof}
We test \eqref{eq:mu-integrated} with $\mu\varphi^2$ and obtain
\begin{align}\label{eq:weak-test}
 nH\int_\Sigma\mu^3\varphi^2\,dV
 &\le \int_\Sigma |A|^2\mu^2\varphi^2\,dV
       -\int_\Sigma\varphi^2|\nabla\mu|^2\,dV\notag\\
 &\hspace{1.5em}
       -2\int_\Sigma\mu\varphi
                  \langle\nabla\mu,\nabla\varphi\rangle\,dV.
\end{align}
On the other hand, we apply \eqref{eq:exterior-stability} to
$\mu\varphi$ and obtain
\[
 \int_\Sigma |A|^2\mu^2\varphi^2\,dV
 \le \int_\Sigma\varphi^2|\nabla\mu|^2\,dV
 +2\int_\Sigma\mu\varphi
                  \langle\nabla\mu,\nabla\varphi\rangle\,dV
 +\int_\Sigma\mu^2|\nabla\varphi|^2\,dV.
\]
Hence we get
\eqref{eq:stability-cancellation}.

For the second estimate, we denote the left hand side by $I=\int_\Sigma\mu^3\chi^p\,dV$ and take
$\varphi=\chi^{p/2}$. \eqref{eq:stability-cancellation} gives
\begin{align*}
     nH I &\le\frac{p^2}{4} \int_\Sigma\mu^2\chi^{p-2}|\nabla\chi|^2\,dV\\
     & =\frac{p^2}{4} \int_\Sigma (\mu^3\chi^p)^{1-2/p}
   \bigl(\mu^{3-p/2}|\nabla\chi|^p\bigr)^{2/p}\,dV\\
   &\le\frac{p^2}{4}I^{1-2/p}\left(\int_\Sigma\mu^{3-p/2}|\nabla\chi|^p\,dV\right)^{2/p}.
\end{align*}
If $I=0$, then the estimate follows. Otherwise, we obtain
\begin{equation*}
 I\le\left(\frac{p^2}{4nH}\right)^{p/2}
       \int_\Sigma\mu^{3-p/2}|\nabla\chi|^p\,dV.
\end{equation*}
Since $p\ge8$ and $\mu\ge H$, we have
\[
 \mu^{3-p/2}
 =\rho\mu^{4-p/2}\le H^{4-p/2}\rho.
\]
Therefore, we get
\eqref{eq:arbitrary-power}.
\end{proof}

\begin{proof}[Proof of Theorem~\ref{thm:main}]
We assume that $H>0$. Since the Morse index is finite, we can choose a compact set $E$ such that
\eqref{eq:exterior-stability} holds outside $E$.
We choose $\zeta\in C^\infty(\Sigma)$ with $0\le\zeta\le1$.
It is zero on a neighborhood of $E$ and one outside a larger compact set $E'$.
We fix a real exponent
\[
 p\ge8,\qquad p>n+1.
\]
We choose $\eta_R$ as the restriction to $\Sigma$ of a smooth ambient cutoff
satisfying
\[
 0\le\eta_R\le1,\qquad
 \eta_R=1\ \text{on }B_R,\qquad
 \eta_R=0\ \text{outside }B_{2R},\qquad
 |\nabla\eta_R|\le C/R.
\]
Properness ensures that $\chi_R=\zeta\eta_R$ has compact support in
$\Sigma\setminus E$. For large $R$, we have $\eta_R=1$ on $E'$ and
$\zeta=1$ on the support of $\nabla\eta_R$. Thus the two terms in
$\nabla\chi_R=\eta_R\nabla\zeta+\zeta\nabla\eta_R$ have disjoint
supports. By Proposition~\ref{prop:weighted-volume}, we obtain
\begin{align}\label{eq:cutoff-error}
 \int_\Sigma\rho|\nabla\chi_R|^p\,dV
 &\le \int_\Sigma\rho|\nabla\zeta|^p\,dV
       +\frac{C^p}{R^p}
        \int_{\Sigma\cap B_{2R}}\rho\,dV\notag\\
 &\le C_0+C_1R^{-p}(2R+H^{-1})^{n+1}
   = C_0+o(1).
\end{align}
Here $C_0<\infty$ because $\nabla\zeta$ has compact support and
$\rho\le H^{-1}$. All constants are independent of $R$.

We apply \eqref{eq:arbitrary-power} to $\chi_R$. Since $\eta_R\to1$
pointwise, we use Fatou's lemma and \eqref{eq:cutoff-error} to obtain
\[
 \int_\Sigma\mu^3\zeta^p\,dV<\infty.
\]
In particular, we obtain
\[
 H^3\operatorname{Area}(\Sigma\setminus E')
 \le\int_{\Sigma\setminus E'}\mu^3\,dV<\infty.
\]
The compact set $E'$ also has finite area. We conclude that
\begin{equation}\label{eq:finite-area}
 \operatorname{Area}(\Sigma)<\infty.
\end{equation}
Since $\Sigma$ is complete and noncompact and has bounded mean
curvature, we apply Frensel's infinite-volume theorem
\cite[Appendix, Theorem~1]{Frensel1996} to obtain
$\operatorname{Vol}_n(\Sigma)=\infty$. This contradicts
\eqref{eq:finite-area}. Hence $H=0$.
\end{proof}

\paragraph{Disclosure of AI tools.}
The author used AI as an interactive aid in developing research ideas, drafting and revising the manuscript, and checking proofs. The author formulated the problem, directed the investigation, and selected the mathematical approaches pursued. The author takes full responsibility for the final content of the paper.

{\small
\let\originalthebibliography\thebibliography
\renewcommand{\thebibliography}[1]{%
 \originalthebibliography{#1}%
 \setlength{\itemsep}{0pt}}
\bibliographystyle{alpha}
\bibliography{main}
}

\medskip

\begingroup
\small
\normalfont
\noindent
Department of Mathematics, University of Rochester\\
Rochester, New York\\
Email address:
\href{mailto:gniu3@ur.rochester.edu}
     {\texttt{gniu3@ur.rochester.edu}}
\par
\endgroup
\end{document}